\documentclass[final,1p,times]{elsarticle}

\usepackage{amssymb}
\usepackage{amsthm}
\usepackage{amscd}
\usepackage{bm}
\usepackage{amsmath}
\usepackage{amsfonts}
\usepackage{graphicx}
\usepackage{color}
\biboptions{sort&compress}
\usepackage[breaklinks,colorlinks]{hyperref}
\numberwithin{equation}{section}
\newtheorem{theorem}{Theorem}[section]

\newtheorem{lemma}[theorem]{Lemma}

\usepackage{mathrsfs}
\usepackage{titletoc}
\usepackage{verbatim}
\usepackage{framed}

\renewcommand{\le}{\left}
\newcommand{\ri}{\right}

\makeatletter

\newcommand{\Rmnum}[1]{\expandafter\@slowromancap\romannumeral #1@}
\makeatother
\usepackage{bbding}

\journal{***}
\begin{document}

\begin{frontmatter}

\title{Fractional Laplace operator and  related Schr\"{o}dinger equations\\ on  locally finite graphs}
\author[author1]{Mengjie Zhang}
\ead{zhangmengjie@mail.tsinghua.edu.cn}
\author[author2]{Yong Lin}
\ead{yonglin@mail.tsinghua.edu.cn}
\author[author3]{Yunyan Yang{\footnote{corresponding author}}}
\ead{yunyanyang@ruc.edu.cn}
\address[author1]{Department of Mathematical Sciences, Tsinghua University, Beijing 100084, China}	
\address[author2]{Yau Mathematical Sciences Center, Tsinghua University, Beijing 100084, China}
\address[author3]{School of Mathematics, Renmin University of China, Beijing, 100872, China}

\begin{abstract}
 In this paper, we first define a discrete version of the fractional Laplace operator  $(-\Delta)^{s}$ through
 the heat semigroup on a stochastically complete,  connected,  locally finite graph. 
 Moreover, we introduce a fractional Sobolev space, 
which is necessary when we study problems involving $(-\Delta)^{s}$. Thirdly, we define the fractional divergence,  and then give another form of $(-\Delta)^s$, which leads to a formula  of integration by parts.
Finally, using the mountain-pass theorem and the Nehari manifold, we obtain multiplicity solutions to a discrete fractional Schr\"{o}dinger equation. We caution the readers that though these existence results are well known in the continuous case, the discrete case
is quite different.
\end{abstract}


\end{frontmatter}

\section{Introduction}

Fractional  Laplace operators and related problems have been a classical topic in functional and harmonic analysis,
and they  have very important applications in various fields,
such as water waves   \cite{A32,A42,A35},
optimization   \cite{A37},
crystal dislocation  \cite{A8,A47},
conservation laws   \cite{A9},
singular set of minima of variational functionals   \cite{A69,A55},
stratified materials  \cite{A81},
minimal surfaces   \cite{A20},
gradient potential theory  \cite{A71,T8}
and so on.
We highly recommend some classical books and papers that discuss this theme in detail, such as
 \cite{A1, A, A78,A12}  and the many references   therein.

  For the concept of fractional Laplace operators on the Euclidean space, readers can refer to \cite{Kwa15}, which collects ten equivalent definitions.
Here,  we are solely concentrating on a single definition provided by the heat semigroup.
To be exact,    let $ \mathscr{X} $ be any of the Lebesgue space $ {L}^p$ with $p \in [1,+\infty)$,   the space $ {C}_0$ of continuous functions vanishing at infinity, or  the space $ {C}_{bu}$ of bounded uniformly continuous functions.
For any $s\in(0,1)$,  the fractional Laplace  operator   $(-\Delta)^s:   D_{\mathscr{X}}((-\Delta)^s)  \rightarrow \mathscr{X}  $  is assigned as
\begin{align}\label{delta}
 (-\Delta)^s u(x)=\frac{s}{\Gamma(1-s)} \int_0^{+\infty}
\left(u(x)-e^{t \Delta}u(x)\right) t^{-1-s} d t,
\end{align}
where    $ \Gamma(\cdot)$ is the Gamma function and $e^{t \Delta}$
denotes the heat semigroup of the Laplace operator $\Delta$.

As we know, the fractional  Schr\"{o}dinger equation appears in a series of problems involving molecules, atoms, nuclei, and so on, and the results are very realistic. In particular, in Euclidean space,  the fractional  Schr\"{o}dinger equation was introduced as a result of extending the Feynman path integral by Laskin in \cite{B1, B2}.  In fractional quantum mechanics, such an equation is of particular interest in the study of particles on stochastic fields \cite{B3}.  It is formulated as
\begin{align}\label{e}
 (-\Delta)^s  u+h(x) u=f(x, u)  \ \text{ in } \ \mathbb{R}^N,
\end{align}
where  $s \in(0,1)$,  $h(x): \mathbb{R}^N \rightarrow \mathbb{R}$ is the potential function and
$f(x,y): \mathbb{R}^N \times \mathbb{R} \rightarrow \mathbb{R}$ is the nonlinear term.
 There are a lot of research results in this direction.  Felmer-Quaas-Tan
   \cite{B6} considered \eqref{e} with $h(x)\equiv 1$ and $f(x, u)$ satisfies some suitable assumptions, and obtained its classical positive solutions.
 Dipierro-Palatucci-Valdinoci \cite{B7}  proved that there exists positive  solutions of \eqref{e}, in the case of $h(x)\equiv 1$ and $f(x, u)=|u|^{p-1}u$ with $1 < p <  {2N}/({N-2s})- 1$.
For the same nonlinearity $f$,  Cheng  \cite{B8} showed the existence of ground states for  \eqref{e}, when $h$ satisfies a coercivity assumption, i.e., $h (x)\rightarrow +\infty$ as $|x| \rightarrow +\infty$.
Moreover, under the same coercivity assumption of $h$, Secchi   \cite{B9} extended the results of ground states
to general nonlinearity $f(x, u)$. For further details about fractional Schr\"{o}dinger equations, we refer to \cite{B16,7,3,5} and references therein.

 Note that the study of partial differential equations on graphs or metric spaces has a long history, see for example, a book \cite{Keller-1}
for a nice survey. Recently, in a series of work
\cite{A-Y-Y-1, A-Y-Y-2, A-Y-Y-3}, Grigor’yan-Lin-Yang obtained solutions of various nonlinear equations on graphs by the variational methods.
In particular, in \cite{A-Y-Y-3}, they extended the results for Schr\"odinger equations in the Euclidean space \cite{doO,Adi-Yang,YangJFA} to locally finite graphs, by using the mountain-pass theorem.
 Since then, there has been a lot of work around this interesting topic, see for examples \cite{Zhang-Zhao,Hua-Xu,S-Y-Z-2,Han-Shao-Zhao,Liu-Zhang,De-C}.

 Now we concern existing work for  the discrete form of the   fractional Laplacian \eqref{delta}  for any $s\in(0,1)$.
 {On the integer lattice  graph  $\mathbb{Z}$,
  \eqref{delta}  is  studied in  \cite{CRSTV2018} and     \cite{KN23},  and  an equivalent form is given as
  \begin{align}\label{Wang}
(-\Delta)^s u(x)
=   \sum_{y \in  \mathbb{Z},y\ne x} K_s(x-y) \big(u(x)-u(y)\big)  ,\ \forall u\in L^2(\mathbb{Z}),
\end{align}
where the discrete kernel $K_s$ is given by
  \begin{align*}
  K_s(x-y)= \frac{s}{\Gamma(1-s)} \int_0^{+\infty}e^{t \Delta}\delta_y(x) t^{-1-s} dt ,  \  x\ne y,
\end{align*}
and $\delta_y$   denotes the  Dirac function. 
On $N$-dimensional  lattice graph  $\mathbb{Z}^N$ ($N\geq 1$),    \eqref{delta}  is  studied in \cite{LR2018}.
In addition, on a stochastically complete graph $G$ with standard measure and weight,   \eqref{delta}  is  studied in \cite{Wang-frac}   and an equivalent form is wrote as  \eqref{Wang}.}
 From \cite{Hua-Lin, Wo09},  a weighted graph $G = (V, E, \mu,  {w})$ is called stochastically complete if
 \begin{align}\label{heat-3}
 \sum_{y\in V} p(t, x, y)\mu(y)=1, \ \forall t >0,\  \forall x\in V,
 \end{align}
where $p(t,x,y)=(e^{t\Delta}\delta_y)(x)/\mu(y)$ is the heat kernel on $G$.
  In fact, any finite graph is stochastically complete and $p(t,x,y)$ can be represented by the eigenvalues and eigenfunctions of   $-\Delta$ on any finite graph.
    Starting by this point,  in \cite{Z-L-Y}, we get an explicitly represented of $(-\Delta)^s$ on finite graphs.

 For other forms of the Euclidean fractional Laplace operator, there are also corresponding discrete versions and related problems,  see for examples   \cite{LM2022,  JZ22, XZ20}. {A more interesting version is the definition through the Fourier transform. We refer  readers to  \cite{CRSTV2015, Tarasov, CRSTV2017, LM2022}  and references therein. }
In this direction, L. Wang \cite{WangLiDan}
  considered  ground state solutions  of    \eqref{e} on $\mathbb{Z}^N$, where  the operator $ (-\Delta)^s $ is defined by the Fourier transform in \cite{LM2022},   $h$ is a bounded periodic positive function and $f$ is a periodic function satisfying certain conditions.  \\

 Our aim is twofold.
 One is to extend   \eqref{Wang} to a  more general stochastic complete graph  with arbitrary measure  and weight for any bounded function;
   the other is to establish a fractional Sobolev space $W^{s,2}(V)$, on which we prove multiplicity of solutions to (\ref{e})
    by using the mountain-pass theorem and the Nehari manifold respectively.

\section{Notations and main results}

We  begin by reviewing some basic notions and concepts about graphs.
Let $V$ be a set of vertices, and $E=\{[x, y]: x,\ y\in V  \text{ satisfy }  x\sim y\}$ be a set of edges, where $x\sim y$ means $x$ is connected to $y$ by an edge.
 For any vertex $x\in V$,  assign its measure $\mu:V\rightarrow \mathbb{R}^{+}$ with $x\mapsto \mu(x)$.
 For any edge $[x,\,y]\in E$,  denote its weight $w:E\rightarrow \mathbb{R}^{+}$ with  $[x,\,y]\mapsto  w_{xy}$ and $w_{xy}=w_{yx}$ for any $x\sim y$. Then the quadruple $G = (V, E, \mu, w)$ refers to a weight graph.
  A graph is called connected if any two vertices can be connected via finite edges, and is called locally finite if each vertex has finite neighbors.
   Moreover, let     $C(V)$ be the space of all real-valued functions defined on  $V$, and  $C_c(V)$ be the space of all compactly supported functions.
	The integral of a function $u\in C(V)$ is defined by $
		\int_{V}ud\mu=\sum_{x\in V}u(x)\mu(x).$
For any  $p\in [1,+\infty)$,    ${L^p(V)}$   is a linear  function space   with the finite norm
 $\|u\|_{ p }=\left(\sum_{x\in  V}|u(x)|^p\mu(x)\right)^{ {1}/{p}},$
and  ${L^\infty(V)}$ is a linear space of all bounded functions on $V$ with the finite norm
$  \|u\|_{ \infty }=\sup_{x\in  V}|u(x)|.$

 In this paper, we always assume that  the fractional exponent $s\in(0,1)$ and  the graph $G = (V, E, \mu,  {w})$ is    a stochastically complete,  connected,  locally finite graph.
Define the fractional Laplace  operator $(-\Delta)^{s}$ acting on a function $u\in L^\infty(V)$ as
\begin{align}\label{ee-s}
 (-\Delta)^s u(x)=\frac{s}{\Gamma(1-s)} \int_0^{+\infty}
\left(u(x)-e^{t \Delta}u(x)\right) t^{-1-s} d t,
\end{align}
where the heat semigroup
 \begin{align}\label{e2}
e^{t \Delta}u(x)
=\sum_{y\in V} p(t, x, y) u(y) \mu(y) \end{align}
and the heat kernel $p(t,x,y)$  satisfies \eqref{heat-3}. We next give a  discrete version  of  the fractional Laplace operator.

\begin{theorem}\label{T1}
The fractional Laplace  operator $(-\Delta)^s$ in \eqref{ee-s} can be expressed as
 \begin{align} \label{Delta-s}
(-\Delta)^s u(x) 
   =\frac{1}{\mu(x)}\sum_{y \in  V, \, y \neq  x}W_s(x, y) \left(u(x)-u(y)\right),  \ \forall u\in L^\infty(V),
\end{align}
 where the discrete kernel  $W_s(x, y)$ is a symmetric positive function  defined by
\begin{align}\label{Hs}
W_s(x, y)= \frac{s} {\Gamma(1-s)} \mu(x)\mu(y)\int_0^{+\infty}  p(t, x, y) t^{-1-s} d t,\ \forall x \neq y \in  V.
\end{align}
Moreover,  {there exists a  positive constant $   C_{x,s} $ depending only on $x$ and $s$, such that}
 \begin{align}\label{W}
 { \sum_{y \in  V, \, y \neq  x}W_s(x, y) \leq    C_{x,s} .}
\end{align}
  \end{theorem}
\noindent   {\textbf{Remark.} This key inequality  \eqref{W} ensures   that the fractional Laplacian $(-\Delta)^s$ in \eqref{Delta-s} is well defined for $u\in L^\infty(V)$.}

 Secondly,  we introduce the fractional gradient and {a fractional Sobolev space}.
Note that only the case of countably infinite  $V$ is considered in this paper,
since the case of finite $V$ has been considered in \cite{Z-L-Y}.  Then we  rewrite
\begin{align*}
V=\{ y_1, y_2, \cdots, y_k, \cdots\}.
\end{align*}
 For convenience,
    we define a function $\widetilde{W}_s(x, y): V\times V \rightarrow \mathbb{R}$ 
    by
   \begin{align*}
\widetilde{W}_s(x, y)=   \left\{\begin{aligned}
    &W_s(x, y),  &\text{ if } x \neq y,\\
    &0,&\text{ if } x= y,\\
    \end{aligned}\right.
    \end{align*}
    where $W_s(\cdot, \cdot)$ is defined as in \eqref{Hs} and satisfies \eqref{W}.
Inspired by \cite{S-Y-Z-1},   we assign  the  fractional gradient $\nabla^s$ acting on a function $u\in L^\infty(V)$  as a vector-valued function
\begin{equation}\label{e-7}
  \nabla^s u (x)=\left(\sqrt{\frac{\widetilde{W}_s(x, y_1)}{2\mu(x)}}(u(x)-u(y_1)), \cdots,\sqrt{\frac{\widetilde{W}_s(x, y_k)}{2\mu(x)}}(u(x)-u(y_k)),\cdots\right). 
\end{equation}
Hence the inner product of  the fractional gradient is assigned  as  for any $u$ and $v\in   {L^\infty(V)}$,
\begin{align}\label{grad-s}
 \nabla^su \nabla^sv   (x)
  =   \frac{1}{2\mu(x)}\sum_{y \in  V, \, y \neq x} W_s(x, y) (u(x)-u(y))(v(x)-v(y)),
  \end{align}
and the length of $\nabla^su(x)$    is written  as $|\nabla^su|  (x)=  \sqrt{\nabla^s u  \nabla^s u(x)}$.
{Noting  the property of $W_s(x, y)$ in \eqref{W},    we know that \eqref{grad-s} and $|\nabla^su|  $  are well defined.}
Moreover,  {assuming  $\inf_{x \in V}\mu(x)>0$, }
{ we define  a  fractional Sobolev space} as
 \begin{align*}
W^{s,2}(V) =\left\{u \in  {L^\infty(V)}: \int_{V} \left(|\nabla^su|^2 + u ^2\right)  d\mu
<+\infty\right\},
\end{align*}
which is equipped with   a norm $$
\|u\|_{s,2} = \left(\int_{V}  \left(|\nabla^su|^2 + u ^2\right)  d\mu \right)^{\frac{1}{2}}.$$
It is clear that $W^{s,2}(V)$ is an inner product space  equipped with an inner product
$$\langle u, v\rangle_{s,2} = \int_{ V} \left ( \nabla^s u \nabla^s v+u v  \right )d\mu .$$
Our second result is

\begin{theorem}\label{Sobolev}
The fractional Sobolev space $W^{s,2}(V)$ satisfies the following properties: \\ [0.75ex]
(i) $ C_c(V) \subseteq  W^{s,2}(V)$;\\ [0.75ex]
 (ii) 
 $W^{s,2}(V)$ is a Hilbert space;\\ [0.75ex]
 (iii) 
 $W^{s,2}(V)$ is continuously  embedded in $L^q(V)$ for any $2\leq q\leq+\infty$.
\end{theorem}

 {Thirdly,  we define the fractional divergence,  and then give another form of $(-\Delta)^s$, which leads to a formula  of integration by parts.}
{Set a vector space
\begin{align*}
\mathbb{R}^{\infty}=\left\{\mathbf{x}=\left(x_1, x_2, \cdots, x_i,\cdots\right): x_i \in \mathbb{R}, \ i=1,2, \cdots\right\},
\end{align*}
which is equipped with the coordinate unit vectors $\{ \mathbf{e}_1,  \mathbf{e}_2,\cdots,  \mathbf{e}_i, \cdots\}$. Here   $$ \mathbf{e}_i =(e_{i1},e_{i2},\cdots,e_{ii},\cdots)  \in {\mathbb{R}^{\infty}}$$ with $  {e}_{ii} =1$  and $ {e}_{ij}= 0$ if $j\not=i$.
Assign the set of all $\mathbb{R}^{\infty}$-valued functions as
\begin{align*}
V^{\mathbb{R}^{\infty}}=\left\{\mathbf{u}: V  \rightarrow \mathbb{R}^{\infty}: \mathbf{u} =\left(u_1 , u_2 , \cdots,u_i,\cdots\right),  \ u_i\in C(V),\ i=1,2, \cdots\right \}.
\end{align*}
Define  a vector-valued functions space
 \begin{align*}
\mathscr{L}^2(V)=\left\{ \mathbf{u} \in V^{\mathbb{R}^{\infty}} : \int_V|\mathbf{u}(x)|^2 d \mu<+\infty\right\},
 \end{align*}
where  $|\mathbf{u}(x)|= (\sum_{i=1}^{+\infty} u_i ^2(x) )^{{1}/{2}}$
with $\mathbf{u} =\left(u_1 , u_2 , \cdots,u_i,\cdots\right)\in V^{\mathbb{R}^{\infty}}$.
  Since  $\delta_x\mathbf{e}_i\in\mathscr{L}^2(V)$ for every $x\in V$ and $i\geq 1$,    $ \mathscr{L}^2(V) $  is  nonempty. It is evident that the space
  $\mathscr{L}^2(V)$ is an inner product space
equipped with an inner product
\begin{align*}
\langle\mathbf{u},\mathbf{v}\rangle_{\mathscr{L}^2(V)}= \sum_{ x\in V}  \mu (x)  \, \mathbf{u}(x)  \cdot \mathbf{v}(x)  ,
\end{align*}
where $\mathbf{u}(x)  \cdot \mathbf{v}(x)=\sum_{ i=1}^{+\infty} u_i(x)v_i(x)$ with $\mathbf{u} =\left(u_1 , u_2 , \cdots,u_i,\cdots\right)$ and $\mathbf{v} =\left(v_1 , v_2 , \cdots,v_i,\cdots\right)$ in ${\mathscr{L}^2(V)}$.}
  {Moreover, the fractional divergence ``$\mathrm{div}_s$"  acting on  a vector-valued function   $\mathbf{f} \in { \mathscr{L}^2(V)}$
is assigned as
\begin{align}\label{div}
\int_V (\mathrm{div}_s \mathbf{f})\, \varphi \,d \mu=-\int_V \mathbf{f}\cdot  \nabla ^s\varphi  \,d \mu, \  \forall \varphi \in C_c(V).
\end{align}
For any  $x \in V$,   taking $ \varphi$ as the Dirac function $\delta_x/{\mu(x)}$, we obtain
\begin{align}\label{div-s}
\mathrm{div}_s  \mathbf{f} (x)=-\frac{1}{\mu(x)}\sum_{y\in V} \mu(y) \nabla^s \delta_x  (y)    \cdot \mathbf{f}(y),\  \forall  \mathbf{f} \in { \mathscr{L}^2(V)}.
\end{align}
Our third result contains not only an equivalent form of $(-\Delta)^s$ , but also
 a formula  of integration by parts on  the graph, which is  fundamental when we use methods in calculus of variations.}
 \begin{theorem}\label{repres}
 Suppose $ u\in W^{s,2}(V)$.  Then   the fractional   Laplace  operator  $(-\Delta)^s$  in \eqref{Delta-s} can be expressed as
\begin{align}\label{Def-s}
(-\Delta)^s u =-\mathrm{div}_s  \nabla^s u,\end{align}
where $\nabla^s$ and ${\rm div}_s$ are defined by (\ref{e-7})
 and (\ref{div-s}) respectively.
Moreover, there holds
 \begin{align}\label{part} \int_{V} \varphi (-\Delta)^s u  d\mu =  \int_{V}\nabla^s \varphi  \nabla^s u d\mu, \   \forall   \varphi \in C_c(V).\end{align}
 \end{theorem}

Finally,  
we consider a fractional Schr\"{o}dinger equation
\begin{align}\label{eS1}
(-\Delta)^s u+h(x) u=f(x, u) \ \text{ in }\ V,
\end{align}
where $(-\Delta)^s $ is defined as in \eqref{Delta-s} for any $s \in(0,1)$. 
 The potential function $ h: V \to \mathbb{R} $  satisfies
\begin{align}\label{H-1} h_0=\inf_{x \in V}h(x)>0
\end{align}
 and for any fixed $x_0\in V$,
\begin{align}\label{H-2}
 h (x)\rightarrow +\infty \  \text{ as }\ \mathrm{d}(x,x_0) \rightarrow +\infty ,
\end{align}
 where $\mathrm{d}(x,x_0)$   represents the minimum number of edges required for the path between
 $x$ and $x_0$.   Let  $\mathscr{H}_s$  be a   completion of $C_c(V)$ under the norm
 \begin{align*}
\|u\|_{\mathscr{H}_s}=\left(\int_{V} (|\nabla^s u|^2+h u^2 ) d \mu\right)^{\frac{1}{2}}.
\end{align*}
 Concerning the existence of solutions to  \eqref{eS1}  in the space ${\mathscr{H}_s}$,
we have the following:

\begin{theorem}\label{T3}
Let $G = (V, E, \mu,  {w})$ be a stochastically complete,  connected,  locally finite graph.
Suppose that $\inf_{x \in V}\mu(x)>0$, $h: V  \rightarrow \mathbb{R}$  satisfies \eqref{H-1} and \eqref{H-2},     and  $ f(x, y): V \times \mathbb{R} \to \mathbb{R} $      satisfies $(\mathrm{F}_1)$--$(\mathrm{F}_4)$ as follows:\\[0.75ex]
$\left(\mathrm{F}_1\right)$  $f(x,y)$ is  continuous with respect to $y\in \mathbb{R}$, and $f(x, 0)=0$ for any $x\in V$;\\  [0.75ex]
$\left(\mathrm{F}_2\right)$   for any positive constant $M$, there exists a constant $C_M$ such that
\begin{align*}|f(x, y)| \leq  C_M,\  \forall  (x,y) \in V\times[-M, M];\end{align*}
$\left(\mathrm{F}_3\right)$ there exists a constant  $\alpha>2$  such that
\begin{align*}
0 <\alpha F(x, y)  \leq y f(x, y),\  \forall  (x,y) \in V\times \mathbb{R} \backslash \{0\},
\end{align*}
where $F(x, y)= \int_0^y f(x, t) d t$  is the primitive function of $f$;\\ [0.75ex]
$(\mathrm{F}_4)$    there holds
\begin{align*} \limsup_{y \rightarrow 0} \frac{ f(x, y)}{y}<\lambda_1  =\inf _{u\in\mathscr{H}_s,u\not\equiv 0} \frac{\|u\|_{\mathscr{H}_s}^2}{\|u\|_2^2}
\end{align*}
uniformly in $x \in V$.\\ [0.75ex]
 Then  \eqref{eS1} has at least two solutions in ${\mathscr{H}_s}$: one is strictly positive, and the other is strictly negative.
\end{theorem}

\begin{theorem}\label{T4}
Under the assumptions of Theorem \ref{T3},
if further $f$ satisfies  \\[0.75ex]
$(\mathrm{F}_5)$ for all $x\in V$,  $ { f(x, y)}/{|y|}$  is    strictly increasing  with respect to   $y\in (-\infty,0)$ and $ (0,+\infty)$. \\[0.75ex]
 Then \eqref{eS1} has at least two ground state solutions in ${\mathscr{H}_s}$: one is strictly positive, and the other is strictly negative.
\end{theorem}

Both the proofs of  Theorems \ref{T3}  and  \ref{T4}  are based on the mountain-pass theorem and the Nehari manifold.  We remark that   $f(x,y)$ need not be odd with respect to $y$ as in \cite{WangLiDan}.  For example,
if we set
$$f(x,y)=\left\{ \begin{aligned} &y(e^{y^2}-1), &y\geq 0,\\
&y^3(e^{y^4}-1),& y<0,\end{aligned}\right.
$$
then $f(x,y)$ satisfies the hypotheses $(\mathrm{F}_1)$--$(\mathrm{F}_5)$.\\

 The remaining part of this paper is divided into two sections:
 in Section \ref{S3},  we prove Theorems  \ref{T1}--\ref{repres};
in Section  \ref{S4},   we prove  Theorems \ref{T3}  and  \ref{T4}.
 For simplicity,  we do not distinguish between sequence and subsequence unless necessary. Sometimes, we denote various constants by
 the same letter $C$.

 \section{Proofs of Theorems \ref{T1}--\ref{repres}}\label{S3}

 In this section, we  prove Theorems \ref{T1}--\ref{repres} separately.

\subsection{Proof of Theorem \ref{T1}}

 We first show that the discrete kernel   $W_s(x, y)$ defined by \eqref{Hs}   has a precise upper bound.

\begin{lemma}\label{L1.2}
 For any fixed $x\in V$,
there holds
 \begin{align*}
\sum_{y \in  V, \, y \neq  x}W_s(x, y)\leq  {C_{x,s}} ,
\end{align*}
where the constant
\begin{align}\label{e10}
 {C_{x,s}=\frac{ \mu(x) } {(1-s)\Gamma(1-s)}\cdot\max\left\{ \mu(x)  \max_{t\in   [0,1]  }|\partial_t p(t, x, x) |,\, 1\right\}.}
   \end{align}
\end{lemma}
\begin{proof}
For any {fixed} $x\in V$,  since the vertices set  $V$  is  countably infinite,
then we  rewrite it as
$V=\{ x,  y_2, \cdots, y_k , \cdots\}$.
For any $  t\in [0,+\infty)$, we set
 \begin{align*}
{S_n(t) =\sum_{k=2}^n p(t, x, y_k) \mu(y_k)  \  \text{ and }\
 S(t) =\sum_{k=2}^{+\infty} p(t, x, y_k) \mu(y_k)   .}
\end{align*}
It is clear that
\begin{align}\label{e0}
S_n(t)\rightarrow S(t)\ \text{  point-wise  in  }   t\in [0,+\infty).
\end{align}
From the condition in \eqref{heat-3}, we obtain
  \begin{align}\label{e1}
  S(t) = 1-p(t, x, x) \mu(x)\leq   1,\ \forall t\in [0,+\infty).
\end{align}
If $t\in [0,1]$,  then by using  the  properties of  $p(t, x, y)$ in  \cite{HLLY,Wo09}:
{$p(0,x,y)= \delta_x(y)/{\mu(x)}$ and $p(t, x, y)$ is smooth with respect to $t\in [0,+\infty)$,} we obtain
  \begin{align*}
   S(t)=\left(  p(0, x, x)  -p(t, x, x) \right)\mu(x)
 \leq  C_{x}  t,
  \end{align*}
where the constant  $C_{x} =\max\left\{ \mu(x)  \max_{t\in   [0,1]  }|\partial_t p(t, x, x) |,\, 1\right\}$.
This together with  \eqref{e0} and  Lebesgue's control convergence theorem leads to that
\begin{align}\label{e4}
 \lim_{n\rightarrow +\infty}\int_0^1 t^{-1-s}S_n(t) d t
 = \int_0^1 t^{-1-s}S(t) d t\leq
\frac{ C_{x} }{1-s}.
\end{align}
 Combining \eqref{e1} and \eqref{e4},  we derive that
  \begin{align}\label{e12}
  \nonumber\lim_{n\rightarrow +\infty}\int_0^{+\infty} t^{-1-s}S_n(t) d t=&\int_0^{+\infty}  t^{-1-s}S(t) d t\\
  \nonumber\leq&\frac{ C_{x} }{1-s}+\int_1^{+\infty}  t^{-1-s} d t\\
\leq& \frac{ C_{x}  }{s(1-s)} .
  \end{align}
  Therefore, it follows from  \eqref{Hs} and \eqref{e12}  that
\begin{align*}
\nonumber\sum_{y \in  V, \, y \neq  x}W_s(x, y) &=\lim_{n\rightarrow +\infty} \sum_{k=2}^n W_s(x, y_k) \\
\nonumber&= \frac{s \mu(x)} {\Gamma(1-s)}\lim_{n\rightarrow +\infty}\int_0^ {+\infty} t^{-1-s} S_n(t) d t\\
&\leq    C_{x,s} ,
\end{align*}
 where the constant $ C_{x,s} $ is given by \eqref{e10}.
And then we get the lemma.
\end{proof}
 \textbf{\textit{Proof of Theorem  \ref{T1}}.}
For any fixed $x\in V$,   we  rewrite $V=\{ x,  y_2, \cdots,  y_k , \cdots\}$.
 And for any $u\in L^\infty(V)$ and any  $  t\in [0,+\infty)$,    denote
\begin{align*}
{S_n(t,u) =\sum_{k=2}^n p(t, x, y_k) \mu(y_k) (u(x) - u(y_k) ) \ \text{ and }\
 S(t,u)  =\sum_{k=2}^{+\infty} p(t, x, y_k) \mu(y_k)(u(x) - u(y_k) )  .}
\end{align*}
Similar to the proof in Lemma \ref{L1.2},   we obtain
  \begin{align*}
\lim_{n\rightarrow +\infty}\int_0^{+\infty} t^{-1-s}S_n(t,u) d t
 =\int_0^{+\infty} t^{-1-s}S(t,u) d t.
  \end{align*}
   This together with   \eqref{heat-3},  \eqref{ee-s} and \eqref{e2} leads to
\begin{align*}
(-\Delta)^s u(x)  &=\frac{s}{\Gamma(1-s)} \int_0^{+\infty}\left(u(x)-e^{t \Delta}u(x)\right) t^{-1-s} d t\\
&=\frac{s}{\Gamma(1-s)} \int_0^{+\infty}\left(1\cdot u(x)-\sum_{y \in  V}p(t, x, y) u(y)\mu(y)\right) t^{-1-s} d t\\
 &=\frac{s}{\Gamma(1-s)} \int_0^{+\infty} t^{-1-s}  S (t,u)  d t\\
 &=\frac{s}{\Gamma(1-s)}\lim_{n\rightarrow +\infty}\int_0^{+\infty} t^{-1-s}S _n(t,u) d t\\
 & = \frac{1}{ \mu(x)}\sum_{k=2}^{+\infty} W_s(x, y_k) (u(x) - u(y_k) )  .
 \end{align*}
Then  \eqref{Delta-s} follows.  In view of  Lemma \ref{L1.2}, we conclude this theorem.
$\hfill\Box$\\

 \subsection{Proof of Theorem \ref{Sobolev}} \label{S2}

 In this subsection, we investigate some properties of  $W^{s,2}(V) $.  We divide Theorem \ref{Sobolev} into the following three lemmas and prove each one.
    At first,
 we  point out that   $W^{s,2}(V) $  is not empty.

 \begin{lemma}\label{L1.1}
  There holds $ C_c(V) \subseteq  W^{s,2}(V)$.
  \end{lemma}
  \begin{proof}
Let   $\mathrm{supp} (u)=\{x\in V: u(x)\not=0\}$ be  the support set of   function $u$.
For any $u\in C_c(V)$,  it is easy to know that
$\int_{V}  u ^2  d\mu$ is bounded.
On the other hand, it follows from  \eqref{W}  that
 \begin{align*}
 \sum_{x\in \mathrm{supp} (u)} \sum_{y \in  V\backslash \mathrm{supp} (u) } W_s(x, y) u^2(x)
  \leq  C,
\end{align*}
and then there holds
 \begin{align*}
 \sum_{x\in \mathrm{supp} (u)} \sum_{y \in  V\backslash \mathrm{supp} (u) } W_s(x, y) u^2(x)
&= \sum_{y \in  V\backslash \mathrm{supp} (u) } \sum_{x\in \mathrm{supp} (u)}W_s(x, y) u^2(x)\\
&=  \sum_{x \in  V\backslash \mathrm{supp} (u) } \sum_{y\in \mathrm{supp} (u)}W_s(x, y) u^2(y).
\end{align*}
This together with    \eqref{grad-s} leads to
 \begin{align*}
\int_{V} |\nabla^su|^2    d\mu=&\frac{1}{2}\sum_{x\in V} \sum_{y \in  V, \, y \neq x} W_s(x, y) (u(x)-u(y))^2\\
=&\frac{1}{2}\sum_{x\in \mathrm{supp} (u)} \sum_{y \in  \mathrm{supp} (u), \, y \neq x} W_s(x, y) (u(x)-u(y))^2
 +\frac{1}{2}\sum_{x\in \mathrm{supp} (u)} \sum_{y \in  V\backslash \mathrm{supp} (u) } W_s(x, y) u^2(x) \\
&+\frac{1}{2}\sum_{x\in V\backslash \mathrm{supp} (u)} \sum_{y \in  \mathrm{supp} (u) } W_s(x, y) u^2(y) \\
=&\frac{1}{2}\sum_{x\in \mathrm{supp} (u)} \sum_{y \in  \mathrm{supp} (u), \, y \neq x} W_s(x, y) (u(x)-u(y))^2
 + \sum_{x\in \mathrm{supp} (u)} \sum_{y \in  V\backslash \mathrm{supp} (u) } W_s(x, y) u^2(x).
\end{align*}
Since $u\in C_c(V)$ and \eqref{W},  we obtain that $\int_{V} |\nabla^su|^2    d\mu$ is bounded, and thus $ u\in W^{s,2}(V)$.
\end{proof}

Secondly,  we prove that $W^{s,2}(V)$  is complete.

\begin{lemma}\label{P1}
  $W^{s,2 }(V)$ is a Hilbert space.
\end{lemma}
\begin{proof}
 {It is sufficient to prove the completeness of  $W^{s,2 }(V)$.}
 If $\left\{u_n\right\}$ is a Cauchy sequence in $W^{s,2}(V)$, then $\left\{u_n\right\}$ and    $\left\{\nabla^s u_n\right\}$ are   two Cauchy sequences in  $L^2(V)$ and $\mathscr{L}^2(V)$ respectively.
 On the one hand,  $L^2(V)$  is a Hilbert space  from \cite{S-Y-Z-1}. Then   there exists a function $u\in L^2(V)$ such that $u_n\rightharpoonup u$ weakly in $L^2(V)$ as $n\rightarrow+\infty$. This implies that
\begin{align}\label{7}
u_n\rightarrow u\ \text{  point-wise  in  }   V.
\end{align}
On the other hand,   it is not difficult to know that the space $\mathscr{L}^2(V)$ is a Hilbert space from  the theory of functional analysis.  Then  there exists a vector-valued  function $\mathbf{g}\in  \mathscr{L}^2(V)$ such that $\nabla^s u_n\rightharpoonup \mathbf{g}$ weakly in
$  \mathscr{L}^2(V)$, and $\nabla^s u_n\rightarrow \mathbf{g}$   point-wise  in $V$.
This together with \eqref{e-7}  and \eqref{7} leads to   $$\mathbf{g}(x)=\nabla^su(x),  \ \forall x\in V.$$ Hence the function $u\in W^{s,2}(V)$ and  the  Cauchy sequence $\left\{u_n\right\}$   converges to $u$ in $W^{s,2}(V)$. We conclude that   $W^{s,2 }(V)$ is complete, and thus it is a Hilbert space.
\end{proof}

Thirdly,
we introduce a fractional Sobolev embedding.

\begin{lemma}\label{Le1-1}
$W^{s,2}(V)$ is continuously  embedded in $L^q(V)$ for any $2\leq q\leq+\infty$.
\end{lemma}

\begin{proof}
For any $u\in W^{s,2}(V)$,
it is obvious that $ \|u\|_2\leq \|u\|_{s,2}$, then  $W^{s,2}(V)$ is embedded in $L^2(V)$. Moreover,    it follows from $ \mu_0=\inf_{x \in V}\mu(x) >0$ that
 \begin{align*}
\|u\|^2_\infty
{ \leq  }    \sum_{x\in V} u^2(x)  \leq   \frac{1}{\mu_0} \|u\|^2_2 \leq   \frac{1}{\mu_0}
\|u\|_{s,2}^2,
\end{align*}
and thus $W^{s,2}(V)$ is embedded in $L^{\infty} (V)$.
Considering the two estimates above, we get
for any $2<q<+\infty$,
\begin{align*}
  \|u\|_q^q
= \int_V   |u|^{q-2}  u^2 d\mu\leq \|u\|_\infty^{q-2} \|u\|_2^2 \leq {\mu_0^{\frac{2-q}{2}}}\|u\|_{s,2}^q.
\end{align*}
Hence $W^{s,2}(V)$ is embedded in $L^q(V)$ for any $2\leq q\leq+\infty$.
\end{proof}

 In summary, we prove Theorem \ref{Sobolev} in connection with  Lemmas \ref{L1.1}--\ref{Le1-1}.

\subsection{Proof of Theorem  \ref{repres}}

 On the one hand,  for any fixed $x \in V$,  we  rewrite  $V=\{ x,  y_2,  \cdots, y_k, \cdots\}$. It follows from \eqref{grad-s} that
\begin{align*}
 \nabla^s \delta_x   \nabla^su (x)
  &=   \frac{1}{2\mu(x)}\sum_{y \in  V, \, y \neq x} W_s(x, y) (u(x)-u(y)),\\
 \nabla^s \delta_x   \nabla^s u   (y_k)
 & =  \frac{1}{2\mu(y_k)} W_s( x,y_k) ( u(x)-u(y_k)),\ \forall k=2,3,\cdots.
  \end{align*}
 This together with \eqref{div-s} leads to
\begin{align*}
-\mathrm{div}_s\nabla ^su(x)=&\nabla ^s\delta_x \nabla ^su(x)+\frac{1}{\mu (x)}\sum_{k=2}^{+\infty}{\mu (y_k)}\nabla ^{s}\delta_x  \nabla ^su(y_k)
\\
=&\frac{1}{2\mu(x)}\sum_{y \in  V, \, y \neq x} W_s(x, y) (u(x)-u(y))+\frac{1}{2\mu (x)}\sum_{k=2}^{+\infty}{W_s(x,y_k)(u(x)-u(y_k))}
\\
=&{\frac{1}{\mu (x)}}\sum_{y\in V,\,y\ne x}W_s(x,y)(u(x)-u(y)),
\end{align*}
 and then  \eqref{Def-s} holds from  \eqref{Delta-s}.
 On the other hand, it is clear that\eqref{part} follows from \eqref{div} and  \eqref{Def-s}.  Then Theorem  \ref{repres} is proved.  

\section{Fractional Schr\"{o}dinger equation}\label{S4}

In this section,
we employ the theory of critical point to prove Theorems \ref{T3}  and  \ref{T4}  respectively.  We first do some necessary preliminary analysis. After that, we discuss the existence of strictly positive and strictly negative solutions in two separate subsections.

 \subsection{Compactness analysis}

At first, we give a compactness result related to $\mathscr{H}_s$.

\begin{lemma}\label{Le4-1}
If $\inf_{x \in V}\mu(x)>0$, then
the space $\mathscr{H}_s$ is weakly pre-compact and compactly embedded in  $L^q(V)$ for any  {$2\leq q\leq+\infty$.}
Moreover, if $\{u_n\}$ is bounded in $\mathscr{H}_s$, then   there exists a function $u\in\mathscr{H}_s$ such that
	    \begin{align*}\left\{
        \begin{aligned}
			&u_n\rightharpoonup u \ \text{ weakly in }  \mathscr{H}_s,\\
			&u_n\rightarrow u \ \text{ strongly in } L^q(V),\ \forall  q\in [2,+\infty],\\
            &u_n\rightarrow u\ \text{  point-wise  in  }   V.
        \end{aligned}\right.
      	\end{align*}
    \end{lemma}

\begin{proof}
 It follows from $ h_0>0$ that  {for any $u\in \mathscr{H}_s$,}
 \begin{align*}
\|u\|^2 _2
\leq \frac{1}{h_0} \sum_{x\in V} h(x)u^2(x)\mu(x)
\leq   C \|u\|^2 _{\mathscr{H}_s},
\end{align*}
and then   there holds
 \begin{align*}
\|u\|_{s,2}^2&
=\|u\|_{\mathscr{H}_s}^2+ \int_{V} (1-h)u^2   d\mu\\
&\leq \|u\|_{\mathscr{H}_s}^2+ (1- h_0)\|u\|^2 _2
\\
&\leq   C \|u\|^2 _{\mathscr{H}_s},
\end{align*}
 which leads to $\mathscr{H}_s$ is  embedded in  $W^{s,2}(V)$.
Due to Lemma \ref{Le1-1}, we conclude  that $\mathscr{H}_s$ is  embedded in $L^q(V)$ for any $2\leq q\leq+\infty$.

On the other hand, it is obvious that ${\mathscr{H}_s}$ is a Banach space.
And it follows  from $ h_0>0$ that  ${\mathscr{H}_s}$ is a  continuously embedded  subspace of $W^{s,2}(V)$:
 \begin{align*}
\|u\|_{s,2}^2
=\|u\|_{\mathscr{H}_s}^2+ \int_{V}\left( \frac{1}{h}-1\right)hu^2   d\mu
 \leq 
C \|u\|^2 _{\mathscr{H}_s}.
\end{align*}
 Then we conclude that
   ${\mathscr{H}_s}$ is a reflexive Banach space because its complete norm is induced by an inner product.
Therefore,  for any bounded sequence $\{u_n\}\subseteq\mathscr{H}_s$,  there exists some function $u\in \mathscr{H}_s$ such that $u_n\rightharpoonup u$ weakly in $\mathscr{H}_s$ as $n\rightarrow+\infty$.
Moreover,    we conclude that   $u_n$ point-wise converges to $u$ as $n\rightarrow+\infty$.
We next prove $u_n\rightarrow u$  strongly in $L^q(V)$ for all $2\leq q\leq+\infty$.

Since $\{u_n\}$ is bounded in $\mathscr{H}_s$ and $u\in\mathscr{H}_s$, there exists some constant $C_1$ such that
	\begin{equation}\label{9}
\int_{V}h(u_n-u)^2  d\mu
\leq 2\int_{V}h(u_n^2+u^2)  d\mu\leq
	2(\|u_n\|^2_{\mathscr{H}_s}+\|u\|^2_{\mathscr{H}_s})\leq C_1.
	\end{equation}
 Let $x_0\in V$ be fixed.
 Taking into account the condition \eqref{H-2},  for any $\epsilon>0$, there exists some large enough $r\geq 1$ such that for $x\in {V\backslash B_r(x_0)} $,
	\begin{equation}\label{10-1}
		h(x)\geq\frac{C_1}{\epsilon}.
	\end{equation}
Combining $(\ref{9})$ and $(\ref{10-1})$,   one has
	\begin{align}\label{11-1}
	\int_{V\backslash B_r(x_0)} (u_n-u)^2d\mu
 \leq   \frac{\epsilon}{C_1}
 \int_{V} h  (u_n-u)^2d\mu
 \leq   \epsilon.
	\end{align}	
	Moreover, it follows from ${  B_r(x_0)}$  is a finite vertexes set  that
	$$\lim_{n\rightarrow+\infty}\int_{B_r(x_0)}(u_n-u)^2d\mu=0.$$
	This together with $(\ref{11-1})$ gives
	 that  $u_n\rightarrow u$ strongly in $L^2(V)$.
Hence it is easy to know that 	$$\|u_n-u\|_\infty^2\leq\frac{1}{\mu_0}\int_{V}|u_n-u|^2d\mu\rightarrow0,$$
where $\mu_0=\inf_{x \in V}\mu(x)$.
For any $2<q<+\infty$, there has
	$$\int_{V}|u_n-u|^qd\mu\leq\|u_n-u\|^{q-2}_\infty\|u_n-u\|_2 ^2 \rightarrow0 .$$
In conclusion,  we get $u_n\rightarrow u$  strongly in $L^q(V)$ for all $2\leq q\leq+\infty$.
 \end{proof}

Secondly,
we denote $u^+=\max\{u,0\}$ and $u^-=\min\{u,0\}$ as the positive and negative part of $u$ respectively.
 We give the following lemma.

\begin{lemma}\label{L2.6}
If $u\in \mathscr{H}_s$, then $u^+,u^-,|u|\in \mathscr{H}_s$.
\end{lemma}

\begin{proof}
For any $u\in \mathscr{H}_s$,   there exists a  Cauchy sequence $\left\{u_n\right\} \subseteq C_c(V)$ such that $\left\|u_n-u\right\|_{\mathscr{H}_s}\rightarrow 0$ as $n \rightarrow +\infty$.
Then we derive
\begin{equation}\label{16-1}
\int_V  |\nabla^s (u_n-u)|^2   d \mu
\rightarrow 0 \ \text{ and }\  \int_V   h (u_n-u)^2  d \mu
\rightarrow 0.
\end{equation}
 For any vectors $ \mathbf{a}$, $ \mathbf{b}\in \mathbb{R}^{\infty}$ satisfying  $ |\mathbf{a}|<+\infty$ and $ |\mathbf{b}|<+\infty$, there holds
 \begin{align*}
 \big| | \mathbf{a}|- |\mathbf{b}|\big| \leq |  \mathbf{a}-\mathbf{b} |.
 \end{align*}     
It follows from  the above inequality,
   the H\"{o}lder inequality,  the  Minkowski  inequality  and  Lemma \ref{Le4-1} that
\begin{align*}
\int_V \left||\nabla^s u_n|^2- |\nabla^s u|^2\right| d \mu
   &\leq  \int_V\left(|\nabla^s u_n|+ |\nabla^s u| \right)  \left|  \nabla^s (u_n -    u) \right| d \mu \\
      &\leq \left(\int_V   \left(|\nabla^s u_n|+ |\nabla^s u| \right)^2d \mu\right)^{\frac{1}{2}}\left(\int_V \left| \nabla^s\left( u_n-u \right) \right|^2 d \mu\right)^{\frac{1}{2}}\\
          &\leq \left( \left\| \nabla^s u_n   \right\|_{L^2} +\left\|  \nabla^s u  \right\|_{L^2}\right)\left\| \nabla^s\left( u_n-u \right) \right\|_{L^2}\\
                  &\leq   \left\| \nabla^s( u_n  -u) \right\|_{L^2}^2 +2\left\|   u  \right\|_{\mathscr{H}_s} \left\| \nabla^s\left( u_n-u \right) \right\|_{L^2}.
   \end{align*}
This together with \eqref{16-1} leads to
 \begin{align} \label{16}\lim_{n\rightarrow+\infty}\int_V  |\nabla^s  u_n |^2   d \mu
=\int_V  |\nabla^s  u |^2   d \mu.\end{align}

On the one hand,  it is clear that $u^+\geq 0$, $u^-\leq 0$ and  $u^+u^-=0$.
Therefore,   for any $u\in  \mathscr{H}_s$, one has
 \begin{align}\label{13}
 \nonumber  {\nabla ^su^+\nabla ^su^-}(x)   &  =   \frac{1}{2\mu(x)} \sum_{y \in  V, \, y \neq x} W_s(x, y) \left(u^+(x)-u^+(y)\right)\left(u^-(x)-u^-(y)\right)\\
\nonumber &  =-   \frac{ 1}{2\mu(x)} \sum_{y \in  V, \, y \neq x} W_s(x, y) \left(  u^+(x)u^-(y)+  u^+(y)u^-(x)\right)\\
&\geq0,
\end{align}
which implies that
 $|\nabla^su|  \geq |\nabla^su^+| .$
Then for any fixed $x_0\in V$ and $r\geq1$,
\begin{align*}
\int_{ V\backslash B_r(x_0)}\left |\nabla^s (u_n^+-u^+)\right|^2 d \mu
   &\leq 2\int_{ V\backslash B_r(x_0)}  |\nabla^s u_n^+|^2   d \mu +2\int_{ V\backslash B_r(x_0)}  |\nabla^s u^+|^2  d \mu\\
   &\leq 2\int_{ V\backslash B_r(x_0)}  |\nabla^s u_n|^2 d \mu +2\int_{ V\backslash B_r(x_0)} |\nabla^s u|^2  d \mu\\
     &=4 \int_{ V\backslash B_r(x_0)} |\nabla^s u|^2  d \mu+o_n(1),
   \end{align*}
   where     the last equality follows from    \eqref{16}.
On the other hand, from $u_n$ uniformly  converges to $ u$ in ${B_r(x_0)}$, we derive
$$\int_{B_r(x_0)}\left |\nabla^s (u_n^+-u^+)\right|^2 d \mu=o_n(1).$$
Since $x_0$ and $r$ are arbitrary, we can take a sufficiently large $r$ such that
$$\int_{ V\backslash B_r(x_0)}\left |\nabla^s (u_n^+-u^+)\right|^2 d \mu+\int_{B_r(x_0)}\left |\nabla^s (u_n^+-u^+)\right|^2 d \mu\leq \epsilon+o_n(1),$$
and thus we get
\begin{align*}
\lim_{n\rightarrow+\infty}\int_V  \left|\nabla^s (u_n^+-u^+)\right|^2 d \mu=
0.
   \end{align*}
In the same way, we also conclude
\begin{align*}
\lim_{n\rightarrow+\infty}\int_V  h (u_n^+-u^+ )^2 d \mu=0.
   \end{align*}
As a consequence, there is $u^+\in \mathscr{H}_s$. Moreover, we also get $u^-\in \mathscr{H}_s$ following the above method.
Since $\mathscr{H}_s$ is complete, it follows from $|u|=u^+-u^-$ that   $|u|\in \mathscr{H}_s$.
\end{proof}

Define the weak solution   $u\in\mathscr{H}_s$ of \eqref{eS1}:  there holds
	\begin{align}\label{weak-s}
		\int_{V}\le( \nabla^s u\nabla^s\varphi+hu\varphi\ri) d\mu=\int_{V}f(x,u)\varphi d\mu,\ \forall \varphi\in \mathscr{H}_s .
	\end{align}
 According to a simple argument,  we can know that  if $u\in\mathscr{H}_s$ is a weak solution of \eqref{eS1}, it is also a point-wise solution:
  for any fixed $x_0\in V$, we take a test function $\varphi_0(x)\in C_c(V)$ satisfying  $\varphi_0(x_0)=1$,   and  $\varphi_0(x)=0$ in the rest. In the following text, the discussion will focus on the weak solution of  \eqref{eS1}.

\subsection{Strictly positive solutions}

In this subsection,   we discuss the existence of strictly positive solutions of  \eqref{eS1} in $\mathscr{H}_s$.
We next show the pivotal lemma associated with the weak solution.

\begin{lemma}\label{L3.1}
Assume $ f(x, y) \equiv 0$,  for any $ y\leq 0$.
If  $u \in \mathscr{H}_s$ is a nontrivial weak solution of \eqref{eS1}, then it is also
 a strictly positive  solution.
\end{lemma}

\begin{proof}
Since  $u\in\mathscr{H}_s$ is a weak solution of \eqref{eS1}, it is also a point-wise solution.
 Moreover, it follows from  Lemma \ref{L2.6} that  $u^-\in \mathscr{H}_s$.    Hence by taking   $\varphi=u^-$  in  \eqref{weak-s}, and using  $u=u^++u^-$, we have
\begin{align*}
 \int_{V}\left(|\nabla^s   u^-|^2 +\nabla^s u^+\nabla^s   u^-  \right) d \mu+\int_{V}h (u^-)^2 d \mu=0.
\end{align*}
Therefore, it follows from \eqref{13} and $ h_0=\inf_{x \in V}h(x)>0$ that
$$\int_{V}h (u^-)^2 d \mu=0,$$
and then  $u^- \equiv 0$ in $V$.
That is to say $u(x)\geq0$ for all $x\in V$.

  We  claim  that $u(x)>0$ for all $x\in V$. Otherwise, if there exists $u\left(x_0\right)=0$ for some $x_0\in V$, then
 $
(-\Delta)^s  u(x_0)=0
$
from   $u\in\mathscr{H}_s$ is  a point-wise solution of \eqref{eS1}.
In view of  the definition of  $(-\Delta)^s$ in \eqref{Delta-s}, we obtain
\begin{align*}
\frac{1}{\mu(x_0)}\sum_{y \in  V, \, y \neq  x_0}W_s(x_0, y) u(y) =0,
\end{align*}
which leads to $u\equiv0$ in $V$ from   $u\geq0$. This contradicts  with the premise that $u $ is nontrivial, and ends the proof.
\end{proof}

According to   Lemma \ref{L3.1}, we obtain that strictly positive solutions of  \eqref{eS1} are equivalent
to nontrivial weak  solutions of the following equation
  \begin{align}\label{eS1-1}
(-\Delta)^s u+h(x) u=f(x, u^+) \ \text{ in }\ V.
\end{align}
 Define the energy functional $J_+:\mathscr{H}_s\rightarrow\mathbb{R}$ associated to   \eqref{eS1-1} by
	\begin{align*}
		J_+(u)=\frac{1}{2}\int_{V}\le(|\nabla^s u|^2+hu^2\ri)d\mu-\int_{V}F(x,u^+)d\mu.
	\end{align*}
In particular, the critical point of $J_+$  is a solution of  \eqref{eS1-1}, which is given by  \eqref{part}.
Therefore,  this subsection will focus on the nontrivial critical point of $J_+$.

\subsubsection{Mountain-pass type solutions}\label{subsub2}

 Firstly, we need to describe the geometry profile of $J_+$, and we have the following lemma.
\begin{lemma}\label{L3.2}
There exists a  nonnegative function $u \in \mathscr{H}_s$ such that $J_+(t u) \rightarrow-\infty$ as $t \rightarrow+\infty$.
\end{lemma}

\begin{proof}
Because of
$\left(\mathrm{F}_3\right)$, we obtain that for any $x\in V$ and $ y>0$, there have $ F(x, y)>0$ and
$$ \frac{\partial}{\partial y }  \log \frac{F(x, y)} {y^{\alpha}}=\frac{yf(x,y)-\alpha F(x,y)}{y F(x,y)}   \geq 0.$$
Then we obtain that  for any $x\in V$,
 $     {F(x, y)}/{y^{\alpha}}   $ is  a  nondecreasing  function  with respect to $y>0$.
 For any fixed $x_0\in V$, we take a test function $u_0 \in \mathscr{H}_s$ satisfying  $u_0(x_0)=1$,   and $u_0(x)=0$ in the rest. Hence for any  $ t>0$, there always exists some $ y_0\in(0, t]$ such that
\begin{equation}\label{2}{F(x_0, t)}\geq \frac{F(x_0, y_0)} {y_0^{\alpha}}{t^{\alpha}}.\end{equation}
Therefore, it follows from  \eqref{W} and \eqref{2} that
\begin{align*}
J_+(t u_0)   &= \frac{t^2}{2} \sum_{y \in  V, \, y \neq x_0} W_s(x_0, y)+\frac{t^2}{2}h(x_0) \mu\left(x_0\right)- F\left(x_0, t\right)\mu\left(x_0\right)\\
&\leq   t^2  C_{x_0,s}
 +  t^2  h(x_0) \mu\left(x_0\right)-{t^{\alpha}} \frac{F(x_0, y_0)} {y_0^{\alpha}}\mu (x_0 )\\
&\leq t^2  C_{x_0,s}   \left( 1  +   h(x_0)  - t^{\alpha-2}  \frac{F(x_0, y_0)} {y_0^{\alpha}}\right),
\end{align*}
where the constant $ C_{x_0,s}>\mu (x_0)  $ depends only on $x_0$ and $s$.
Since $\alpha>2$, there has $J_+(t u) \rightarrow-\infty$ as $t \rightarrow+\infty$, and then the proof is complete. %
\end{proof}

\begin{lemma}\label{L3.3-2}
 There exists a positive constant $r$ such that $\inf_{\|u\|_{\mathscr{H}_s}=r}J_+(u) >0$.
 \end{lemma}

 \begin{proof}
The conditions $(\mathrm{F}_3)$ and $(\mathrm{F}_4)$ yield that for any $x\in V$,
 there exist  two sufficiently small constants $\epsilon>0$   and $y_0>0$ such that
\begin{align*}
0< F(x, y)  \leq \frac{y}{\alpha} f(x, y)
  \leq  \frac{\lambda_1-\epsilon}{2} y^2, \ \forall  y \in(0,  y_0).
\end{align*}
Moreover, it is obvious  that $F(x, y) \leq  {y^3}/{y_0^3}  F(x, y)$ for any  $ y  \geq  y_0 $.
Then for any $(x, y) \in V \times \mathbb{R}^+$, there has
\begin{align*}
F(x, y) \leq \frac{\lambda_1-\epsilon}{2} y^2+\frac{1}{y_0^3} y^3 F(x, y),
\end{align*}
which together with the definition  of $\lambda_1$
leads to
\begin{align}\label{14}
 \int_{V} F(x, u^+ )d \mu   
 \leq  \frac{\lambda_1-\epsilon}{2\lambda_1}    {\|u \|^2_{\mathscr{H}_s} } +\frac{1}{y_0^3} \int_{V} (u^+)^3 F(x, u^+ ) d \mu,\ \forall u\in {\mathscr{H}_s}.
\end{align}
By the Sobolev embedding in Lemma \ref{Le4-1}, for any function $u\in {\mathscr{H}_s}$ with $ \|u\|_{\mathscr{H}_s} \leq  1$, there exist two positive constant $C_1$ and $C_2$ such that
$\|u\|_{\infty} \leq  C_1$ and
 $\|u \|_3^3 \leq  C_2\|u \|_{\mathscr{H}_s}^3$.
According to $\left(\mathrm{F}_2\right)$ and $\left(\mathrm{F}_3\right)$, we get
\begin{align*}
\max_{y\in [0,C_1]} F(x, y)\leq \max_{y\in [0,C_1]}\frac{y}{\alpha} f(x, y)\leq  {C_3} ,
\end{align*}
where  $C_3$ is a constant depending only on $C_1$.
This yields
$$
\int_{V} (u^+)^3 F(x, u ^+) d \mu \leq C_3 \|u \|_3 ^3 \leq  C_2C_3 \|u \|_{\mathscr{H}_s}^3.
$$
Then it follows from \eqref{14} that
\begin{align*}
\int_{V} F(x, u^+ )d \mu   \leq \frac{\lambda_1-\epsilon}{2{\lambda_1}}  {\|u\|^2_{\mathscr{H}_s} } +\frac{C_2C_3}{y_0^3}\|u \|_{\mathscr{H}_s}^3  .
\end{align*}
 Hence there holds for any $u\in {\mathscr{H}_s}$ with $\|u\|_{\mathscr{H}_s} \leq  1$,
$$
\begin{aligned}
J_+(u)  
& \geq  \frac{1}{2}\|u\|_{\mathscr{H}_s}^2-\frac{\lambda_1-\epsilon}{2{\lambda_1}}  {\|u \|^2_{\mathscr{H}_s} } -\frac{C_2C_3}{y_0^3}\|u \|_{\mathscr{H}_s}^3\\
&  \geq\frac{\epsilon}{2 \lambda_1}\left(1-\frac{2 \lambda_1 {C_2C_3}}{\epsilon{y_0^3}} \|u\|_{\mathscr{H}_s}\right)\|u\|_{\mathscr{H}_s}^2 .
\end{aligned}
$$
Setting $r=\min \{1, \epsilon y_0^3 /\left(4 \lambda_1 C_2C_3\right) \}$, we obtain $J_+(u) 
>0$ for all $u\in {\mathscr{H}_s}$ with $\|u\|_{\mathscr{H}_s}=r$. This completes the proof.
 \end{proof}

We now prove that $J_+$ satisfies the  Palais-Smale condition.

\begin{lemma} \label{L3.4-2}
For any $c \in \mathbb{R}$,  the functional $J_+$ satisfies the $(\mathrm{PS})_c$   condition. Namely, if $\{u_n\} \subseteq \mathscr{H}_s$ satisfies $J_+\left(u_n\right) \rightarrow c$ and $J_+'\left(u_n\right) \rightarrow 0$ as $n \rightarrow+\infty$, then there exists  a function $u_0 \in \mathscr{H}_s$ such that  $u_n \rightarrow u_0$ in $\mathscr{H}_s$.
\end{lemma}
\begin{proof}
On the one hand,
since $\{u_n\} \subseteq \mathscr{H}_s$ satisfies $J_+\left(u_n\right) = c+o_n(1)$, there holds
$$\left\|u_n\right\|_{\mathscr{H}_s}^2  =2 \int_{V} F(x, u_n ^+) d \mu+2 c+o_n(1) ,  $$
which together with $\left(\mathrm{F}_3\right)$ leads to
\begin{align}\label{15}
\left\|u_n\right\|_{\mathscr{H}_s}^2
 \leq  \frac{2}{\alpha} \int_{V} f(x, u_n ^+) u_n ^+ d \mu+2 c+o_n(1).
\end{align}
On the other hand, from $J_+'\left(u_n\right) = o_n(1)$, we obtain
  for any $\varphi \in \mathscr{H}_s $,
\begin{align}\label{1}
\int_{V} f(x, u_n^+ ) \varphi d \mu= \int_{V} \left( \nabla^s u_n \nabla^s \varphi+h u_n \varphi\right) d \mu +o_n(1) .
\end{align}
By taking the test function $\varphi=u_n$ in \eqref{1}, one gets
$
\int_{V} f(x, u_n ^+) u_n  ^+d \mu=\left\|u_n\right\|_{\mathscr{H}_s}^2+o_n(1),
$
which together with \eqref{15} leads to
$$
\left\|u_n\right\|_{\mathscr{H}_s}^2
  \leq  \frac{2}{\alpha}\left\|u_n\right\|_{\mathscr{H}_s}^2+2 c+o_n(1) .
$$
Then $u_n$ is bounded in $\mathscr{H}_s$ from $\alpha>2$.

 From  Lemma \ref{Le4-1},   $u_n$ is   bounded on $V$, and
  there exists a function $u_0\in\mathscr{H}_s$ such that $\|u_n-u_0\|_2=o_n(1)$ and
 \begin{align}\label{3}
\int_{V} \left( \nabla^s u_0 \nabla^s  (u_n-u_0)+h u_0  (u_n-u_0)\right) d \mu =o_n(1) .\end{align}
Moreover, in view of $\left(\mathrm{F}_2\right)$ and $(\mathrm{F}_4)$, we split the integral as follows:
\begin{align*}
\left|\int_{V} f(x, u_n ^+)\left(u_n-u_0\right) d \mu\right| &\leq
\int_{|u_n|\leq y_0} \left|f(x, u_n )\right|\left|u_n-u_0\right|  d \mu + \int_{|u_n|> y_0} \left|f(x, u_n )\right|\left|u_n-u_0\right| d \mu\\
&\leq
\int_{| u_n |\leq y_0} \lambda_1 |u_n|  \left|u_n-u_0\right|  d \mu
+C  \int_{|u_n |> y_0}\frac{|u_n|}{y_0}  \left|u_n-u_0\right| d \mu\\
& \leq \left(\lambda_1+\frac{C}{y_0}\right)
\|u_n \|_2 \|u_n-u_0\|_2=o_n(1),
\end{align*}
where $y_0>0$ is a sufficiently small constant. Then by taking the test function $\varphi=u_n-u_0$ in \eqref{1}, we get
\begin{align}\label{4}
\int_{V} \left( \nabla^s u_n \nabla^s  (u_n-u_0)+h u_n  (u_n-u_0)\right) d \mu =\int_{V} f\left(x, u_n^+\right)\left(u_n-u_0\right) d \mu+o_n(1)
=o_n(1) .
\end{align}
Hence it follows from \eqref{3} and \eqref{4} that
$\left\|u_n-u_0\right\|_{\mathscr{H}_s}
=o_n(1),$ and thus the lemma   follows.
\end{proof}

From the three lemmas above,  we check that
the functional $J_+$ satisfies the geometric conditions of the mountain-pass theorem in \cite{A-R}.

\begin{lemma}\label{Lm}
If   $f(x,y)$   satisfies
$\left(\mathrm{F}_1\right)$--$\left(\mathrm{F}_4\right)$, then \eqref{eS1} has a strictly positive solution in ${\mathscr{H}_s}$.
\end{lemma}
\begin{proof}
From Lemmas \ref{L3.2} and \ref{L3.3-2}, we know that there exists $u^*\in \mathscr{H}_s$ with  $\|u^*\|_{\mathscr{H}_s}>r>0$ such that $$\inf_{\|u\|_{\mathscr{H}_s}=r}J_+(u) >J_+(0)=0>J_+(u^*) $$ for some positive constant $r$.
These together with Lemma \ref{L3.4-2} lead to that $J_+$ satisfies all assumptions of the
mountain-pass theorem. 
 Denote
 $$
\Gamma=\Big\{\gamma \in C([0,1], \mathscr{H}_s): \gamma(0)=0, \gamma(1)=u^*\Big\}.
$$
Using the mountain-pass theorem, we conclude that
$$
c=\inf _{\gamma \in \Gamma} \max _{t\in [0,1] } J_+(\gamma(t))
$$
is a critical value of $J_+$.
In particular, there exists some $u \in \mathscr{H}_s$ such that $J_+(u)=c$. 
 Since
$
J_+(u)=c >0,
$
we have  $u \not \equiv 0$.
Then recalling Lemma \ref{L3.1}, we conclude this lemma.
\end{proof}

\subsubsection{Ground state solutions}\label{subsub3}

  Write the Nehari manifold as $\mathscr{N}_+=\left\{u \in \mathscr{H}_s\backslash\{0\}:
 \langle J_+'(u), u\rangle =0 \, \right\}$, namely
	\begin{equation}\label{N}
		\mathscr{N}_+=\left\{u \in \mathscr{H}_s\backslash\{0\}:   \int_{V} (|\nabla^s u|^2+hu^2 ) d \mu=\int_{V}f(x,u^+ )u^+  d\mu \right\}.
	\end{equation}
Moreover, we define
	\begin{align*}
		c_+=\inf_{u\in\mathscr{N}_+}J_+(u).
	\end{align*}
 If $c_+$ can be achieved by some function $u\in\mathscr{N}_+$, then $u$ has the least energy among all functions belonging to the Nehari manifold and in fact, $u$ is a ground state solution of \eqref{eS1-1}.

 Firstly,   we are ready to prove an important lemma as follows.

\begin{lemma}\label{L4}
If  $  u \in \mathscr{H}_s$ with $u^+ \not\equiv 0$,
then  there exists a unique $t_0 >0$ such that
$t_0 u \in \mathscr{N}_+$ and $$J_+(t_0 u)=\max_{t>0}J_+(t u).$$   Moreover, if $u\in \mathscr{N}_+$, then $J_+( u)=\max_{t>0}J_+(t u)$.
\end{lemma}

\begin{proof}
For any $  u \in \mathscr{H}_s$ with $u^+ \not\equiv 0$ and any $t>0$, we denote
 \begin{align*}
g(t):=J_+(t u)=\frac{t^2}{2}\|u\|_{\mathscr{H}_s}^2-\int_{V}F(x,tu^+ )d\mu,\end{align*}
and assign  the derivative
$g'(t) 
 = t   \bar{g}(t),$
where $$ \bar{g}(t)=\|u\|_{\mathscr{H}_s}^2 - \int_{\{x\in V:\,u^+\neq 0\}} (u^+)  ^2 \frac{f(x, tu ^+) }{tu ^+} d\mu.$$
 By assumption $(\mathrm{F}_5)$, we conclude that $\bar{g}(t)$ is strictly 	decreasing  with respect to $t\in(0,\,+\infty)$.
On the one hand,
similar to the proof of Lemma  \ref{L3.2},   we obtain that  for any $x\in V$,
 $     {F(x, y)}/{y^{\alpha}}   $ is  a  nondecreasing  function  with respect to $y>0$.
 Hence for any $y_0>0$,  there is
  \begin{align*}
  F(x, tu^+) \geq \frac{  F(x, y_0 u^+)}{  y_0  ^\alpha}    t ^\alpha, \    \forall t>y_0.
  \end{align*}
Then it follows from $\left(\mathrm{F}_3\right)$ that there exists a constant  $\alpha>2$  such that
\begin{align*}
  (u^+)^2 \frac{ f(x,  t  u^+) }{  (tu^+) } \geq     \frac{\alpha F(x,  t  u^+)}{  t   ^2}  \geq
   \frac{  \alpha t ^{\alpha -2}}{ y_0^\alpha  }    F(x, y_0 u^+) .
\end{align*}
As a consequence,   for any $t>y_0$, there is
 \begin{align*}
 \bar{g}(t)
  \leq  \|u\|_{\mathscr{H}_s}^2-   \frac{    \alpha }{ y_0^\alpha  }  t ^{\alpha -2}   \int_V  F(x, y_0 u^+) d\mu.
\end{align*}
 Noting $\alpha>2$ and $   \int_V     F(x, y_0 u^+) d\mu>0$, we  derive
$ \bar{g}(t) \rightarrow-\infty$ as $t \rightarrow+\infty$.
 On the other hand,
 from
$(\mathrm{F}_4)$, there has
$$ \limsup _{t \rightarrow 0^+ } \frac{ f(x, tu ^+)}{tu^+ }<\lambda_1 
\leq \frac{\|u\|_{\mathscr{H}_s}^2 }{\|u\|_2^2}.$$
Taking into account the above estimate, it is easy to know that
\begin{align*}
\liminf _{t \rightarrow 0^+ }\bar{g}(t)  =&  \|u\|_{\mathscr{H}_s}^2-\int_{V}(u^+ )^2\limsup _{t \rightarrow 0^+ } \frac{f(x, tu^+ ) }{tu^+ } d\mu \\
 >& \|u\|_{\mathscr{H}_s}^2- \lambda_1\|u\|^2_2   \\
\geq &0.
\end{align*}
Hence there exists a unique $t_0 \in(0, +\infty)$ such that $\bar{g}(t_0)=0$.
Therefore,  we conclude that   $t_0 $ is the unique root for $g'(t)=0$,    and then $g(t_0)=\max_{t>0}g(t)$ and $t_0 u \in \mathscr{N}_+$.

 Moreover, if $u\in \mathscr{N}_+$, it is obvious that $  u \in \mathscr{H}_s$ with $u^+ \not\equiv 0$. By an easy calculation, there has
\begin{align*}
g'(t)&= t  \|u\|_{\mathscr{H}_s}^2-\int_{V}f(x, tu^+  ) u^+   d\mu,
\end{align*}
which implies $g'(1)=0$ from $u\in \mathscr{N}_+$.  Due to the above work, we conclude
 $g(1)=J_+(u)=\max_{t>0}J_+(tu)$, and thus we get the lemma.
\end{proof}


Secondly, we use the above tools to prove that $c_+$ can be achieved in $\mathscr{N}_+$.

\begin{lemma}\label{L3}
There exists some $u_s \in \mathscr{N}_+$ such that $J_+ \left( u_s \right)=c_+>0.$
\end{lemma}
\begin{proof}
 Taking a sequence $\left\{u_n\right\} \subseteq \mathscr{N}_+$ such that $J_+ \left(u_n\right)\rightarrow c_+$
   as ${n \rightarrow +\infty}$, we obtain
$$\left\|u_n\right\|_{\mathscr{H}_s}^2  =2 \int_{V} F(x, u_n ^+) d \mu+2 c_++o_n(1) ,  $$
which together with $\left(\mathrm{F}_3\right)$ and $ \|u_n\|_{\mathscr{H}_s}^2=\int_{V}f(x, u_n^+ )u_n ^+d\mu$ leads to
\begin{align*}
\left\|u_n\right\|_{\mathscr{H}_s}^2
 \leq
 \frac{2}{\alpha} \left\|u_n\right\|_{\mathscr{H}_s}^2+2 c_++o_n(1).
\end{align*}
 Since $\alpha>2$,
this gives that
   $\left\{u_n\right\}$ is bounded in ${\mathscr{H}_s}$.
From Lemma \ref{Le4-1},   there exists some $u_s\in {\mathscr{H}_s}$ such that
	    \begin{equation}\label{18}
\left\{
        \begin{aligned}
			&u_n\rightharpoonup u_s \ \text{ weakly in }  \mathscr{H}_s,\\
			&u_n\rightarrow u_s \ \text{ strongly in } L^q(V),\ \forall  q\in [2,+\infty],\\
            &u_n\rightarrow u_s\ \text{  point-wise  in  }  V.
        \end{aligned}\right.
      	\end{equation}
 By   $(\mathrm{F}_4)$,
 there are  two sufficiently small constants $\epsilon>0$   and $y_0>0$ such that
\begin{align} \label{e24}
 f(x, y)  \leq  (\lambda_1-\epsilon) y, \quad  \forall  (x,y)\in V\times[0,y_0).
\end{align}
 We can derive that  $\|u_n^+\| _\infty\geq y_0>0$,   and then  $u_s^+ \not\equiv 0$.

Set
\begin{align*}
\Omega=\{x\in V:  0\leq u_n ^+(x)< y_0 \  \text{ and } \ 0\leq u_s^+ (x)< y_0   \}.
\end{align*}
On the one hand,
 we get that
\begin{align*}
 \left|\int_{u_s^+}^ {u_n ^+} f(x, t) d t\right|
 \leq (\lambda_1-\epsilon)  \left|\int_{u_s^+}^ {u_n ^+}   t d t\right|
 \leq  \frac{\lambda_1-\epsilon}{2}     ( u_s^++u_n ^+ ) \left|u_s -u_n \right|,\ \forall x\in \Omega
\end{align*}
from  \eqref{e24} and  the inequality   $|a^+ -b^+ |\leq |a-b|$ for any $a,\, b \in  \mathbb{R}$.
This together with  the H\"{o}lder inequality  implies
\begin{align}\label{19-1}
\nonumber\left|\int_{\Omega}\left(F(x, u_n^+ )-F(x, u_s^+)\right) d \mu\right|&\leq  \int_{\Omega}\left| F(x, u_n^+ )-F(x, u_s^+)  \right|d \mu
 \\
 \nonumber&\leq  \int_{\Omega} \left|\int_{u_s^+}^ {u_n ^+} f(x, t) d t\right|   d \mu
 \\
 &\leq \frac{\lambda_1-\epsilon}{2}  \left( \|u_s\|_2+\|u_n \|_2 \right)
   \left\|u_s -u_n \right\|_2.
\end{align}
On the other hand,  $\left\{u_n\right\}$ is bounded in $L^\infty(V)$   by some positive constant $M>y_0$ from   \eqref{18}.
And we know that
\begin{align*}
V\backslash\Omega=\{x\in V:  y_0\leq  u_n ^+(x)\leq M    \}\cup \{x\in V: y_0\leq u_s^+ (x)\leq M   \}.
\end{align*}
 Then from $\left(\mathrm{F}_2\right)$ and the H\"{o}lder inequality, we get
\begin{align*}
\nonumber\left|\int_{V \backslash\Omega}\left(F(x, u_n^+ )-F(x, u_s^+)\right) d \mu\right|
\nonumber &\leq C_M   \int_{V \backslash\Omega}|{u_n^+} -{u_s^+}|d \mu  \\
&\leq \frac{ C_M }{y_0} \left(\int_{u_n ^+\geq y_0} |u_n^+||{u_n ^+} -{u_s^+ }|d \mu +\int_{u_s^+\geq y_0} |u_s^+|  |{u_n^+ } -{u_s^+}|d \mu\right) \\
&\leq \frac{ C_M }{y_0}\left( \|u_n\|_2+\|u_s \|_2 \right)
 \|{u_n } -{u_s }\|_2.
\end{align*}
This together with \eqref{18} and  \eqref{19-1} leads to
\begin{align}\label{20}
\lim _{n \rightarrow +\infty} \int_{V} F(x,u_n^+ ) d \mu=\int_{V} F(x,u_s^+) d \mu .
\end{align}
 Since $u_n\rightharpoonup u_s$  weakly in  $\mathscr{H}_s$
and $u_n \in \mathscr{N}_+$, we get
\begin{equation}\label{17}
\left\|u_s\right\|_{\mathscr{H}_s}^2 \leq  \lim _{n \rightarrow +\infty}  \left\|u_n\right\|_{\mathscr{H}_s}^2=\lim  _{n \rightarrow +\infty}\int_{V}f(x, u_n^+ )u_n ^+ d\mu=\int_{V}f(x,u_s^+ )u_s  ^+d\mu.
\end{equation}
According to \eqref{17} and \eqref{20}, we derive
\begin{align}\label{21}
J_+(u_s) \nonumber =&\frac{1}{2}\|u_s\|_{\mathscr{H}_s}^2-\int_{V}F(x,u_s^+ )d\mu \\
\nonumber  \leq&  \lim_{n \rightarrow +\infty} \left(\frac{1}{2}\|u_n\|_{\mathscr{H}_s}^2-\int_{V}F(x,u_n ^+)d\mu\right) \\
\nonumber =&\lim _{n \rightarrow +\infty} J_+\left(u_n\right) \\
 =&c_+.
\end{align}

We   claim that $u_s \in \mathscr{N}_+$.  Suppose not,  we have  $\left\|u_s\right\|_{\mathscr{H}_s}^2<\int_{V}f(x,u_s^+ )u_s ^+d\mu$ and $u_s^+ \not\equiv 0$. By  $u_s^+ \not\equiv 0$ and  Lemma \ref{L4},  there exists a unique $t_0 >0$ such that 
$t_0 u_s \in \mathscr{N}_+$. Then there has $J_+(t_0 u_s) \geq c_+$.
Combining \eqref{20} and \eqref{17},  we have
\begin{align*}
J_+(t_0 u_s)  & =\frac{t_0^2}{2}\|u_s\|_{\mathscr{H}_s}^2-\int_{V}F(x,t_0u_s ^+)d\mu\\
&< \frac{t_0^2}{2}\int_{V}f(x,u_s ^+)u_s^+ d\mu-\int_{V}F(x,t_0u_s^+ )d\mu\\
&=\lim _{n \rightarrow +\infty} \left( \frac{t_0^2}{2}\int_{V}f(x,u_n^+ )u_n^+ d\mu-\int_{V}F(x,t_0u_n ^+)d\mu\right)\\
&=\lim _{n \rightarrow +\infty} \left( \frac{1}{2}\|t_0u_n\|_{\mathscr{H}_s}^2-\int_{V}F(x,t_0u_n ^+)d\mu\right)\\
&=\lim _{n \rightarrow +\infty} J_+(t_0 u_n).
 \end{align*}
Since $u_n \in \mathscr{N}_+$, we also have $J_+(u_n)=\max _{t >0} J_+(tu_n)$ from Lemma \ref{L4}.  Hence there holds
$$c_+\leq J_+(t_0 u_s)< \lim _{n \rightarrow +\infty} J_+ \left(u_n\right)=c_+,$$
which is a contradiction, and thus $u_s \in \mathscr{N}_+$.
 Therefore,   we conclude  $J_+ \left( u_s \right)=c_+$ from \eqref{21}.

In addition, since $u_s \in \mathscr{N}_+$ and \eqref{N}, there are $  \|u_s\|_{\mathscr{H}_s}\neq 0$ and
$
  \|u_s\|_{\mathscr{H}_s}^2=\int_{V}f(x,u_s ^+)u_s^+ d\mu.
$
As a consequence, it follows from   $\left(\mathrm{F}_3\right)$  that
$$J_+(u_s)  
\geq\frac{1}{2}\|u_s\|_{\mathscr{H}_s}^2-\int_{V}f(x,u_s^+ )\frac{u_s ^+}{\alpha}d\mu  =\left(\frac{1}{2}-\frac{1}{\alpha} \right)\|u_s\|_{\mathscr{H}_s}^2.$$
Since $\alpha>2$,
this gives   $c_+=\inf_{u\in\mathscr{N}_+}J_+(u)>0$, and then the proof is complete.
\end{proof}

Because of Lemmas
 \ref{L4} and \ref{L3}, we conclude that there exists some $u_s \in \mathscr{N}_+$ such that
\begin{align}\label{22-s}
\max_{t>0}J_+(t u_s)=J_+(u_s)=c_+ >0.
\end{align}
We next prove that $u_s$ is the critical point of $J_+$.
\begin{lemma}\label{L2}
There holds $J_+'(u_s)=0$.
\end{lemma}

\begin{proof}
Inspired by  ( \cite{{S-Y-Z-2}}, Lemma 3.5),  we use a modified argument from (\cite{Adimurthi}, Lemma 3.5) to prove this lemma.
Suppose not,    there exists some $\varphi_s \in  \mathscr{H}_s\backslash \{0\}$ such that $ \langle J_+'(u_s), \varphi_s  \rangle\neq0$.
Without loss of generality, we assume  $ \langle J_+'(u_s), \varphi_s  \rangle<0$.   For any $t$, $\tau \in \mathbb{R}$, we first define a function $g (t,\tau)$ by $$g (t, \tau)=t u_s+\tau \varphi_s  .$$
 By an easy calculation,  there has
$$
\lim _{(t,\tau) \rightarrow(1,0)} \frac{\partial}{\partial \tau} J_+(g (t,\tau))=\left\langle J_+'(u_s), \varphi_s \right\rangle<0.
$$
Hence there exist two sufficiently small constants $0<\epsilon_1<1$ and $\epsilon_2>0$ such that   $ J_+ (g (t,\tau))$ is strictly decreasing with respect to $\tau\in \left[-\epsilon_2, \epsilon_2\right]$   for any fixed $t\in \left[1-\epsilon_1, 1+\epsilon_1\right] $. This together with \eqref{22-s} leads to
\begin{align}\label{23}
J_+(g (t,\tau))<J_+(g (t,0))= J_+(t u_s) \leq c_+,\ \forall \tau \in \left(0, \epsilon_2\right].
\end{align}
     Moreover, for any $t,\tau \in \mathbb{R}$, define
$$
g^*(t,\tau)=\|g  \|_{\mathscr{H}_s}^2-\int_{V} f\left(x, g ^+   \right) g^+     d \mu .
$$
Then there holds $
g^*(t,0)   
  =t^2\bar{g}^*(t) ,$
where $$
\bar{g}^*(t) = \|u_s\|_{\mathscr{H}_s}^2-\int_{\{x\in V:\,u_s^+\neq 0\}}\frac{ f\left(x, t u_s ^+\right) }{t u_s^+ } ( u_s ^+)^2 d \mu .$$
It follows from $(\mathrm{F}_5)$ that $\bar{g}^*(t) $ is strictly decreasing with respect to $t$.
This together with the fact of $u_s \in \mathscr{N}_+$ implies
$$
\bar{g}^* (1+\epsilon_1 )<\bar{g}^*(1) =0< \bar{g}^* (1-\epsilon_1 ),$$
namely,
$$
 g^* (1+\epsilon_1 ,0)  < 0<  g ^* (1-\epsilon_1 ,0)  .$$
Applying   the  continuity  of  $g^*(t,\tau)$  with respect to $\tau \in\mathbb{R}$, we obtain that
 there exists a constant $\epsilon_3\in (0,\epsilon_2 )$ such that
 $$ {g^* (1+\epsilon_1 ,\tau)} < 0< { g ^* (1-\epsilon_1 ,\tau)} ,\ \forall \tau \in\left(0, \epsilon_3\right].$$
 Since $g^*(t,\tau)$ is also  continuous  with respect to $t \in\mathbb{R}$,
there exists some $t_{\tau} \in\left(1-\epsilon_1, 1+\epsilon_1\right)$  such that $g^*\left(t_{\tau}, \tau\right)=0$ for all $\tau \in\left(0, \epsilon_3\right].$ As a consequence, we conclude $g \left(t_{\tau}, \tau\right)  \in\mathscr{N}_+$.
Hence it follows from  \eqref{23} that
$$
c_+ \leq J_+\left(g \left(t_{\tau}, \tau\right)\right)< c_+,
$$
which is a contradiction.
\end{proof}

\begin{lemma}\label{Ln}
If $f(x,y)$   satisfies
$\left(\mathrm{F}_1\right)$--$\left(\mathrm{F}_5\right)$, then
\eqref{eS1}   has a strictly positive ground state solution in ${\mathscr{H}_s}$.
\end{lemma}
\begin{proof}
According to Lemmas
 \ref{L4}--\ref{L2},  we conclude that $c_+$ can be achieved by some $u_s \in \mathscr{N}_+$ with $J_+'(u_s)=0$.
Then $u_s$  is not only a ground state solution, but also a nontrivial weak solution of \eqref{eS1-1}.
From Lemma \ref{L3.1}, we conclude this result.
\end{proof}

\subsection{Proofs of Theorems \ref{T3} and \ref{T4}}

In this subsection,  we focus on the existence of strictly negative solutions of  \eqref{eS1} in $\mathscr{H}_s$, and then prove Theorems \ref{T3} and \ref{T4} respectively.
After minor modifications to the proofs of strictly positive solutions, we can also get the existence of strictly negative solutions.
However, for the convenience of readers and the completeness of the article, we give the differences in this subsection.

 \begin{lemma}\label{BL3.1}
 Suppose  $f(x, y) \equiv 0, \ \ \forall  y\geq 0.$
If  $u \in \mathscr{H}_s$ is a nontrivial weak solution of \eqref{eS1}, then it is also
 a strictly negative solution.
\end{lemma}

\begin{proof}
Similar to the proof of Lemma \ref{L3.1},
by taking   $\varphi=u^+$  in  \eqref{weak-s},   we conclude that  $u^+ \equiv 0$ in $V$, i.e.,
  $u(x)\leq0$ for all $x\in V$.  Using reductio ad absurdum, we can obtain   $u(x)<0$ for all $x\in V$, and end the proof.
\end{proof}

From Lemma \ref{BL3.1},  we obtain that strictly negative solutions of  \eqref{eS1} are equivalent
to nontrivial weak  solutions of the following equation
 \begin{eqnarray}\label{BeS1-1}
(-\Delta)^s u+h(x) u=f(x, u^-) \ \text{ in }\ V.
\end{eqnarray}
Moreover, we define the energy functional $J_-:\mathscr{H}_s\rightarrow\mathbb{R}$ associated to   \eqref{BeS1-1} by
	\begin{eqnarray*}
		J_-(u)=\frac{1}{2}\int_{V}\le(|\nabla^s u|^2+hu^2\ri)d\mu-\int_{V}F(x,u^-)d\mu.
	\end{eqnarray*}
In particular, if $u$ is a critical point of $J_-$, then $u$ is a solution of  \eqref{BeS1-1} from  \eqref{part}.
Therefore,  we will focus on the nontrivial critical point of $J_-$ in this subsection.\\


In the next three lemmas, we check that
the functional $J_-$ satisfies the geometric conditions of the mountain-pass theorem.

\begin{lemma}\label{BL3.2}
There exists a nonpositive function $u \in \mathscr{H}_s$ such that $J_-(t u) \rightarrow-\infty$ as $t \rightarrow +\infty$.
\end{lemma}

\begin{proof}
For any $x\in V$ and $ y<0$,
it follows from
$\left(\mathrm{F}_3\right)$  that  there have  $f(x,y)<0$, $ F(x, y)>0$ and
$$ \frac{\partial}{\partial y} \log \frac{F(x, y)} {|y|^{\alpha}}=\frac{|y|f(x,y)+\alpha F(x,y)}{|y| F(x,y)}   \leq 0.$$
Therefore,   for any $x\in V$,
 $     {F(x, y)}/{|y|^{\alpha}}   $ is a   decreasing  function   with respect to $y<0$.
 For any fixed $x_0\in V$, we take a test function $u_0(x)\in \mathscr{H}_s$ satisfying  $u_0(x_0)=-1$,   and $u_0(x)=0$ in the rest. Hence  for any fixed $y_0<0$ and all $t\geq |y_0|$,
\begin{equation}\label{B2}
{F(x_0,- t)}\geq \frac{F(x_0, y_0)} {|y_0|^{\alpha}}{t^{\alpha}}.
\end{equation}
Then from  \eqref{W}  and  \eqref{B2}, we obtain
\begin{align*}
J_-(t u_0)  
\leq t^2 C_{x_0,s}  \left( 1  +   h(x_0)  -t^{\alpha-2}  \frac{\mu(x_0)}{C_{x_0,s}} \frac{F(x_0, y_0)} {|y_0|^{\alpha}}\right),
\end{align*}
where  the constant $C_{x_0,s} >0$ depends only on $x_0$ and $s$.
Since $\alpha>2$, then $J_-(t u) \rightarrow-\infty$ as $t \rightarrow +\infty$.
\end{proof}

\begin{lemma}\label{BL3.3-2}
 There exists a positive constant $r$ such that $\inf_{\|u\|_{\mathscr{H}_s}=r}J_-(u) >0$.
 \end{lemma}

 \begin{proof}
In a similar way of \eqref{14},
  for any $u\in {\mathscr{H}_s}$ with $ \|u\|_{\mathscr{H}_s} \leq  1$, we get there exist  two sufficiently small constants $\epsilon>0$   and $y_0>0$ such that
\begin{align}\label{B14}
\int_{V} F(x, u^- )d \mu
 \leq \frac{\lambda_1-\epsilon}{2\lambda_1}    {\|u \|^2_{\mathscr{H}_s} } +\frac{1}{y_0^3} \int_{V} |u^-|^3 F(x, u^- ) d \mu  .
\end{align}
Moreover, there are two positive constants, $C_1$ and $C_2$, such that
$\|u\|_{\infty} \leq  C_1$ and
 $\|u\|_3^3 \leq  C_2\|u\|_{\mathscr{H}_s}^3$ from  Lemma \ref{Le4-1}. Then it follows from
  $\left(\mathrm{F}_2\right)$ and $\left(\mathrm{F}_3\right)$ that
$$\max_{y\in [-C_1,0]} F(x, y)\leq \max_{y\in [-C_1,0]}\frac{|y|}{\alpha}| f(x, y)|\leq  {C_3} ,$$
where  $C_3$ is a constant depending only on $C_1$.
This yields
$$
\int_{V} |u^-|^3 F(x, u ^-) d \mu \leq C_3 \|u\|_3 ^3 \leq  C_2C_3 \|u\|_{\mathscr{H}_s}^3.
$$
Then it follows from \eqref{B14} that
\begin{eqnarray*}
\int_{V} F(x, u^- )d \mu   \leq \frac{\lambda_1-\epsilon}{2{\lambda_1}}  {\|u \|^2_{\mathscr{H}_s} } +\frac{C_2C_3}{y_0^3}\|u\|_{\mathscr{H}_s}^3  .
\end{eqnarray*}
 Hence we obtain
$$
\begin{aligned}
J_-(u)
  \geq\frac{\epsilon}{2 \lambda_1}\left(1-\frac{2 \lambda_1 {C_2C_3}}{\epsilon{y_0^3}} \|u\|_{\mathscr{H}_s}\right)\|u\|_{\mathscr{H}_s}^2 .
\end{aligned}
$$
Setting $r=\min \{1, \epsilon y_0^3 /\left(4 \lambda_1 C_2C_3\right) \}$, we obtain $J_-(u) 
>0$ for all $u\in {\mathscr{H}_s}$ with $\|u\|_{\mathscr{H}_s}=r$.
 \end{proof}

\begin{lemma} \label{BL3.4-2}
For any $c \in \mathbb{R}$,  the functional $J_-$ satisfies the $(\mathrm{PS})_c$   condition. 
\end{lemma}
\begin{proof}
This proof is similar to Lemma \ref{L3.4-2},
we omit the details and leave the proof to interested readers.
\end{proof}

 \textbf{\textit{Proof of Theorem  \ref{T3}}.}
Lemmas \ref{BL3.2}--\ref{BL3.4-2} lead to that $J_-$ satisfies all   assumptions of the
mountain-pass theorem in \cite{A-R}.
Using the mountain-pass theorem, we conclude that $J_-$
has a nontrivial critical point.
Then recalling Lemma \ref{BL3.1},  we conclude that \eqref{eS1} has a strictly negative solution in ${\mathscr{H}_s}$. This together with the result  in Lemma  \ref{Lm}   leads to this theorem.
 $\hfill\Box$\\


In the next,  we  use the method of Nehari manifold
to prove Theorem   \ref{T4}.
  Write the Nehari manifold as $\mathscr{N}_-=\left\{u \in \mathscr{H}_s\backslash\{0\}:
 \langle J_-'(u), u\rangle =0 \, \right\}$, namely
	\begin{equation*}
		\mathscr{N}_-=\left\{u \in \mathscr{H}_s\backslash\{0\}:   \int_{V} (|\nabla^s u|^2+hu^2 ) d \mu=\int_{V}f(x,u^- )u^-  d\mu \right\}.
	\end{equation*}
It is obvious that for any $u\in\mathscr{N}_-$, there   are  $u^- \not\equiv 0$   and  $\|u\|_{\mathscr{H}_{s}}>0$.
Moreover, we define
	\begin{align*}
		c_-=\inf_{u\in\mathscr{N}_-}J_-(u).
	\end{align*}
 If $c_-$ can be achieved by some function $u\in\mathscr{N}_-$, then $u$ has the least energy among all functions belonging to the Nehari manifold and in fact, $u$ is a ground state solution of \eqref{BeS1-1}.

\begin{lemma}\label{BL4}
If  $  u \in \mathscr{H}_s$ with $u^- \not\equiv 0$,
then  there exists a unique $t_0 >0$ such that
$t_0 u \in \mathscr{N}_-$ and $$J_-(t_0 u)=\max_{t>0}J_-(t u).$$   Moreover, if $u\in \mathscr{N}_-$, then $J_-( u)=\max_{t>0}J_-(t u)$.
\end{lemma}

\begin{proof}
For any $  u \in \mathscr{H}_s$ with $u^- \not\equiv 0$ and any $t>0$, we denote
$
g(t) =J_-(t u)$, 
and assign  the derivative
$g'(t) 
 = t   \bar{g}(t),$
where
\begin{align*}
 \bar{g}(t)
 =\|u\|_{\mathscr{H}_s}^2 - \int_{\{x\in V:\,u^-\neq 0\}} (u^-)  ^2 \frac{f(x, tu ^-) }{tu ^-} d\mu.
\end{align*}
Since $tu ^-<0$  strictly decreases as $t>0$ increases,  and   $ f(x, y)/y$  strictly decreases as $y<0$ increases from $(\mathrm{F}_5)$, then $ {f(x, tu ^-) }/({tu ^-})$  strictly increases  as $t>0$ increases. Hence we conclude that $\bar{g}(t)$ is strictly 	decreasing  with respect to $t\in(0,\,+\infty)$.

On the one hand,
 similar to the proof of Lemma  \ref{BL3.2},   we obtain that  for any $x\in V$,  ${F(x, y)}/{|y|^{\alpha}}   $ is a   decreasing  function   with respect to $y<0$.
 Hence for any $y_0>0$,  we obtain $  F(x, tu^-) \geq  y_0  ^{-\alpha} t ^\alpha  F(x, y_0 u^-)       $   if   $ t>y_0$.
Then it follows from $\left(\mathrm{F}_3\right)$ that there exists a constant  $\alpha>2$  such that
\begin{align*}
  (u^-)^2 \frac{ f(x,  t  u^-) }{  (tu^-) } \geq     \frac{\alpha F(x,  t  u^-)}{  t   ^2}  \geq
   \frac{  \alpha t ^{\alpha -2}}{ y_0^\alpha  }    F(x, y_0 u^-)  , \   \forall t>y_0.
\end{align*}
As a consequence, there has $ \bar{g}(t)  \leq  \|u\|_{\mathscr{H}_s}^2-    t ^{\alpha -2}   \int_V    \alpha { y_0^{-\alpha } } F(x, y_0 u^-) d\mu$ for any $t>y_0$.
 Noting $\alpha>2$ and $   \int_V      \alpha { y_0^{-\alpha } }F(x, y_0 u^-) d\mu>0$,
 we obtain
$\bar{g} (t)\rightarrow-\infty$ as $t \rightarrow+\infty$.
 On the other hand,
 from
$(\mathrm{F}_4)$, there has
$$ \limsup _{t \rightarrow 0^+  } \frac{ f(x, tu ^-)}{tu^- }<\lambda_1
\leq \frac{\|u\|_{\mathscr{H}_s}^2 }{\|u\|_2^2}.$$
Taking into account the above estimate, it is easy to know that
\begin{align*}
\liminf _{t \rightarrow 0 ^+ }\bar{g}(t)  =   \|u\|_{\mathscr{H}_s}^2-\int_{V}(u^- )^2\limsup _{t \rightarrow 0^+ } \frac{f(x, tu^- ) }{tu^- } d\mu
 > 0.
\end{align*}
Hence there exists a unique $t_0 \in(0, +\infty)$ such that $\bar{g}(t_0)=0$.
Therefore,  we conclude that   $t_0 $ is the unique root for $g'(t)=0$,    and then $g(t_0)=\max_{t>0}g(t)$ and $t_0 u \in \mathscr{N}_-$.
 Moreover, if $u\in \mathscr{N}_-$, then $  u \in \mathscr{H}_s$ with $u^- \not\equiv 0$. It follows from an easy calculation that $g'(t)   =t  \|u\|_{\mathscr{H}_s}^2-\int_{V}f( x,tu^-  ) u^-   d\mu $.
This  implies $g'(1)=0$,
and then
 $g(1)=J_-(u)=\max_{t>0}J_-(tu)$.
\end{proof}

 \textbf{\textit{Proof of Theorem  \ref{T4}}.}
Similar to the proof of Lemma  \ref{L3}, we prove that $c_-$ can be achieved in $\mathscr{N}_-$. This together with   Lemma
 \ref{BL4} leads to that there exists  $u_s \in \mathscr{N}_-$ such that
\begin{align*}\max_{t>0}J_-(t u_s)=J_-(u_s)=c_- >0.\end{align*}
Moreover, it is not   difficult to prove that $u_s$ is a critical point of $J_-$ using a similar method of
Lemma \ref{L2}.
Then $u_s$  is not only a ground state solution, but also a nontrivial weak solution of \eqref{BeS1-1}.
From Lemma \ref{BL3.1}, \eqref{eS1} has a strictly negative ground state solution in ${\mathscr{H}_s}$. This together with the result in Lemma \ref{Ln}   proves this theorem.
 $\hfill\Box$\\


\end{document}